\documentclass[titlepage,draft,12pt]{article} 
\usepackage{amsfonts,latexsym,pstricks,pst-plot} 
\usepackage[a4paper]{geometry}
\newcommand{\G}{\mathcal{G}}
\newcommand{\trebar}[1]{\|#1\|}
\newcommand{\lk}{\lambda_{k}}

\newcommand{\re}{{\mathbb{R}}}
\newcommand{\n}{{\mathbb{N}}}
\newcommand{\m}[1]{m(\au{#1}^{2})}
\newcommand{\au}[1]{|A^{1/2}#1|}
\newcommand{\qed}{{\penalty 10000\mbox{$\quad\Box$}}\medskip}

\newtheorem{thm}{Theorem}[section]
\newtheorem{thmbibl}{Theorem}
\newtheorem{rmk}[thm]{Remark}
\newtheorem{prop}[thm]{Proposition}

\newtheorem{lemma}[thm]{Lemma}

\title{Kirchhoff equations from quasi-analytic to spectral-gap data}
\author{Marina Ghisi\medskip\\ 
{\normalsize Universit\`a degli Studi di Pisa} \\
{\normalsize Dipartimento di Matematica ``Leonida Tonelli''}\\
{\normalsize PISA (Italy)}\\
{\normalsize e-mail: \texttt{ghisi@dm.unipi.it}}
\and 
Massimo Gobbino\medskip\\ 
{\normalsize Universit\`a degli Studi di Pisa} \\
{\normalsize Dipartimento di Matematica Applicata ``Ulisse Dini''}\\ 
{\normalsize PISA (Italy)}\\
{\normalsize e-mail: \texttt{m.gobbino@dma.unipi.it}}}
\date{}

\begin{document}
\maketitle
\begin{abstract}
	In a celebrated paper (Tokyo J.\ Math.\ 1984) K.\ Nishihara proved
	global existence for Kirchhoff equations in a special class of
	initial data which lies in between analytic functions and Gevrey
	spaces.  This class was defined in terms of Fourier components
	with weights satisfying suitable \emph{convexity} and
	\emph{integrability} conditions.
	
	In this paper we extend this result by removing the convexity
	constraint, and by replacing Nishihara's integrability condition
	with the simpler integrability condition which appears in the usual
	characterization of quasi-analytic functions.
	
	After the convexity assumptions have been removed, the
	resulting theory reveals unexpected connections with some
	recent global existence results for spectral-gap data.
	
\vspace{1cm}

\noindent{\bf Mathematics Subject Classification 2000 (MSC2000):}
35L70, 35L90.

\vspace{1cm} 

\noindent{\bf Key words:} integro-differential hyperbolic equation,
global existence, Kirchhoff equations, quasi-analytic functions,
spectral-gap data.
\end{abstract}
 
\section{Introduction}

Let $H$ be a separable real Hilbert space.  For every $x$ and $y$ in
$H$, $|x|$ denotes the norm of $x$, and $\langle x,y\rangle$ denotes
the scalar product of $x$ and $y$.  Let $A$ be a self-adjoint linear
operator on $H$ with dense domain $D(A)$.  We assume that $A$ is
nonnegative, namely $\langle Ax,x\rangle\geq 0$ for every $x\in D(A)$,
so that the power $A^{\alpha}x$ is defined provided that $\alpha\geq
0$ and $x$ lies in a suitable domain $D(A^{\alpha})$.

Given $m:[0,+\infty)\to[0,+\infty)$, we consider 
the Cauchy problem
\begin{equation}
	u''(t)+\m{u(t)}Au(t)=0 
	\hspace{2em}\forall t\in[0,T),
	\label{pbm:h-eq}
\end{equation}
\begin{equation}
	u(0)=u_0,\hspace{3em}u'(0)=u_1.
	\label{pbm:h-data}
\end{equation}
 
It is well known that (\ref{pbm:h-eq}), (\ref{pbm:h-data}) is the
abstract setting of the Cauchy-boundary value problem for the
quasilinear hyperbolic integro-differential partial differential
equation
\begin{equation}
	u_{tt}(x,t)-
	m{\left(\int_{\Omega}\left|\nabla u(x,t)\right|^2\,dx\right)}
	\Delta u(x,t)=0
	\hspace{2em}
	\forall(x,t)\in\Omega\times[0,T),
	\label{eq:k}
\end{equation}
where $\Omega\subseteq\re^{n}$ is an open set, and $\nabla u$ and
$\Delta u$ denote the gradient and the Laplacian of $u$ with respect
to the space variables. 

Throughout this paper we assume that equations (\ref{pbm:h-eq}) and
(\ref{eq:k}) are strictly hyperbolic, namely
\begin{equation}
	\mu:=\inf_{\sigma\geq 0}m(\sigma)>0.
	\label{hp:ndg}
\end{equation}

We also assume that $m$ is locally Lipschitz continuous. We never
assume that the operator is coercive or that its inverse is compact.

We refer to the survey \cite{gg:london} and to the references quoted
therein for more details on this equation and its history.  Here we
just recall that, under our assumptions on $m(\sigma)$, problem
(\ref{pbm:h-eq}), (\ref{pbm:h-data}) has a local solution for all initial
data $(u_{0},u_{1})\in D(A^{3/4})\times D(A^{1/4})$.

Existence of global solutions is for sure the main open problem in 
the theory of Kirchhoff equations. A positive answer has been given 
in five different special cases.
\begin{list}{}{\leftmargin 4em \labelwidth 4em}

	\item[(GE 1)] Special forms of the nonlinearity $m(\sigma)$
	(see~\cite{poho-m}).

	\item[(GE 2)] Dispersive equations (see~\cite{das,gh}, and the
	more recent papers~\cite{matsuyama,yamazaki}).

	\item[(GE 3)]  Spectral-gap initial data 
	(see~\cite{manfrin1,manfrin2,hirosawa2,gg:global}).

	\item[(GE 4)] Analytic initial data (see~\cite{as,das-an-2}).  In
	this case the result is actually stronger since it is enough to
	assume that $m(\sigma)$ is continuous and nonnegative.

	\item[(GE 5)] ``Quasi-analytic'' initial data
	(see~\cite{nishihara}).
\end{list}

We refer to the quoted papers for more details on each approach.

In this paper we pursue the path (GE 5), introduced by K.\ Nishihara
in~\cite{nishihara}.  In that paper he proved global existence for
initial data in suitable spaces, defined by imposing the convergence
of some series where the Fourier components of data are multiplied by
weights satisfying suitable \emph{convexity} and \emph{integrability}
conditions (see Section~\ref{sec:nishihara} for the details).

As remarked in~\cite{nishihara}, the weights defining Gevrey spaces 
never fulfill these assumptions, but there are examples of weights 
satisfying these conditions for which the resulting space contains 
non-analytic functions. In other words, Nishihara's spaces are 
expected to be something in between analytic functions and Gevrey 
spaces, and for this reason this result is often referred to as a 
global existence result for quasi-analytic data.

In this paper we prove a similar result \emph{without the convexity
assumption}, and we show that Nishihara's quite strange integrability
condition can be replaced with the more standard integrability condition
which appears in the usual characterization of quasi-analytic
functions.

From the technical point of view, the removal of the convexity
condition requires a new proof of the key estimate
(Proposition~\ref{lemma:est}), which now can no more be established by
means of Jensen type inequalities as in~\cite{nishihara}.

From the point of view of global existence results, apart from
providing a cleaner statement for the beauty of the art, the removal
of the convexity assumption has a somewhat unexpected impact.  Using
some weird weights (which of course do not satisfy the previous
convexity assumptions) we can indeed exhibit examples of spaces
containing functions with low Sobolev regularity where the Kirchhoff
equation is well posed (see Theorem~\ref{thm:sg-data}).  This
reminded us of the spectral-gap global solutions as defined
in~\cite{manfrin1,manfrin2} and then in~\cite{gg:global}.  The
phenomenology is quite similar, but the context and the proof are
completely different, and it doesn't seem so easy to deduce exactly
Theorem~\ref{thm:sg-data} from the known results on spectral-gap
solutions.

Using the same weird weights, we can also show that there are special
unbounded operators for which Kirchhoff equation is well posed in
Sobolev-type spaces such as $D(A^{\alpha+1/2})\times D(A^{\alpha})$
with $\alpha>1/4$ (see Theorem~\ref{thm:sg-op}).  This yields a new
item for the list of global existence results:
\begin{list}{}{\leftmargin 4em \labelwidth 4em}
	\item[(GE 6)]  Special operators.
\end{list}

More important, this adds a new difficulty in the search of a
counterexample to the ``big problem'', namely the global existence in
$C^{\infty}$ in the concrete case, or in $D(A^{\infty})$ for the
abstract equation.  Now we know indeed that any such counterexample
needs to exploit some property of the operator in order to rule out
the special operators to which Theorem~\ref{thm:sg-op} applies.

This paper is organized as follows.  In Section~\ref{sec:statements}
we give rigorous definitions of the functional spaces we need, we
recall Nishihara's work, and we state our main result.  In
Section~\ref{sec:proofs} we prove our main result.  In
Section~\ref{sec:spectral-gap} we explore the connections with
spectral-gap solutions, and we exhibit some strange consequences of
this theory.

\setcounter{equation}{0}
\section{Statements}\label{sec:statements}

\subsection{Functional spaces}\label{sec:functional}

For the sake of simplicity we assume that $H$ admits a countable
complete orthonormal system $\{e_{k}\}$ made by eigenvectors of $A$.
We denote the corresponding eigenvalues by $\lambda_{k}^{2}$ (with
$\lk\geq 0$), so that $Ae_{k}=\lambda_{k}^{2}e_{k}$ for every
$k\in\n$.  Every $u\in H$ can be written in a unique way in the form
$u=\sum_{k=0}^{\infty}u_{k}e_{k}$, where $u_{k}=\langle
u,e_{k}\rangle$ are the components of $u$.  In other words, every
$u\in H$ can be identified with the sequence $\{u_{k}\}$ of its
components, and under this identification the operator $A$ acts
component-wise by multiplication.

We stress that this simplifying assumption is by no means restrictive.
Indeed the spectral theorem for self-adjoint unbounded operators on a
separable Hilbert space (see \cite[Chapter~VIII]{reed}) states that
any such operator is unitary equivalent to a multiplication operator
on some $L^{2}$ space.  More precisely, for every $H$ and $A$ there
exist a measure space $(M,\mu)$, a function $a(\xi)\in L^{2}(M,\mu)$,
and a unitary operator $H\to L^{2}(M,\mu)$ which associates to every
$u\in H$ a function $f(\xi)\in L^{2}(M,\mu)$ in such a way that $Au$
corresponds to the product $a(\xi)f(\xi)$.
		
As a consequence, all the definitions we give in terms of $u_{k}$ and
$\lambda_{k}$ can be extended to the general case by replacing the
sequence of components $\{u_{k}\}$ of $u$ with the function $f(\xi)$
corresponding to $u$, the sequence $\{\lambda_{k}\}$ of eigenvalues of
$A$ with the function $a(\xi)$, and summations over $k$ with integrals
over $M$ in the variable $\xi$ with respect to the measure $\mu$.
Similarly, there is no loss of generality in using components in the
proof of the a priori estimate needed for our existence result
(Theorem~\ref{thm:main}).

Coming back to functional spaces, using components we have that 
$$D(A^{\alpha}):=\left\{u\in H:\sum_{k=0}^{\infty}
\lambda_{k}^{4\alpha}u_{k}^{2}<+\infty\right\}.$$

Let now $\varphi:[0,+\infty)\to[0,+\infty)$ be any function.  Then for
every $\alpha\geq 0$ and $r>0$ one can set
$$\trebar{u}_{\varphi,r,\alpha}^{2}:=\sum_{k=0}^{\infty}\lambda_{k}^{4\alpha}
	u_{k}^{2} \exp\left(\strut r\varphi(\lambda_{k})\right),$$
and then define the \emph{generalized Gevrey spaces} as
$$\G_{\varphi,r,\alpha}(A):= \left\{u\in
H:\trebar{u}_{\varphi,r,\alpha}<+\infty\right\}.$$

These spaces can also be seen as the domain of the operator
$A^{\alpha}\exp\left((r/2)\varphi(A^{1/2})\right)$.  They are Hilbert spaces
with norm $(|u|^{2}+\trebar{u}_{\varphi,r,\alpha}^{2})^{1/2}$, and
they form a scale of Hilbert spaces with respect to the parameter $r$.
They are a natural generalization of the usual spaces of Sobolev,
Gevrey or analytic functions, corresponding to
$\varphi(\sigma)=\log(1+\sigma)$, $\varphi(\sigma)=\sigma^{1/s}$
($s>1$), and $\varphi(\sigma)=\sigma$, respectively.
In~\cite{gg:local} and~\cite{gg:london} it is shown that these spaces
represent the right setting for Kirchhoff equations.

Spaces of \emph{quasi-analytic functions} fit in this framework.
They correspond to weights $\varphi(\sigma)$ which are continuous,
strictly increasing, and satisfy
\begin{equation}
	\int_{1}^{+\infty}\frac{\varphi(\sigma)}{\sigma^{2}}\,d\sigma
	=+\infty.
	\label{defn:qa}
\end{equation}

To be overpedantic, also the continuity and \emph{strict} monotonicity
assumptions on $\varphi$ are not really needed.  Indeed, for every
nondecreasing function satisfying (\ref{defn:qa}), one can always find
a smaller function which is continuous, strictly increasing, and still
satisfies (\ref{defn:qa}).

We refer to~\cite{cc} for more details on quasi-analytic functions in 
the concrete case.

\subsection{Nishihara's work}\label{sec:nishihara}

The following is the main result of \cite{nishihara}, restated and
somewhat simplified using the notations we have just introduced.

\begin{thmbibl}\label{thm:nishihara}
	Let $H$ be a separable Hilbert space, and let $A$ be a nonnegative
	self-adjoint (unbounded) operator on $H$ with dense domain.  Let
	$m:[0,+\infty)\to(0,+\infty)$ be a locally Lipschitz continuous 
	function satisfying the nondegeneracy condition (\ref{hp:ndg}).
	
	Let $\varphi:[0,+\infty)\to[0,+\infty)$ be a continuous and
	strictly increasing function such that, setting 
	$M(\sigma):=e^{\varphi(\sigma)}$, we have that
	\begin{enumerate}
		\item[$(\varphi 1)$] the function 
		$\sigma\to M(\sqrt{\sigma})$ is convex,  
	
		\item  [$(\varphi 2)$] if $M^{-1}(\sigma)$ denotes the inverse function of 
		$M(\sigma)$, then we have that
		$$\int_{1}^{+\infty}\frac{1}{\sigma M^{-1}(\sigma)}\,d\sigma
		=+\infty.$$
	\end{enumerate}
	Let us finally assume that 
	\begin{equation}
		(u_{0},u_{1})\in\G_{\varphi,1,1/2}(A)\times
			\G_{\varphi,1,0}(A).
		\label{hp:nishihara}
	\end{equation}
	
	Then problem (\ref{pbm:h-eq}), (\ref{pbm:h-data}) admits a unique
	global solution
	\begin{equation}
		u\in C^{1}\left([0,+\infty);\G_{\varphi,1,1/2}(A)\right)\cap
		C^{0}\left([0,+\infty);\G_{\varphi,1,0}(A)\right).
		\label{th:reg-sol}
	\end{equation}
\end{thmbibl}

To be more precise, the original statement involved further
assumptions on $M(0)$, and on the spectrum and the inverse of $A$,
which however can be easily removed using arguments that nowadays are
quite standard.

We point out that, in contrast with other results for Kirchhoff
equations, this solution lies in a fixed Hilbert \emph{space} instead
of a Hilbert \emph{scale} (namely in (\ref{th:reg-sol}) the radius
$r=1$ is the same for all times).

We also remark that in general it is not possible to replace $r=1$ in
(\ref{hp:nishihara}) with a smaller value of $r$.  The point is that,
when we replace $\varphi(\sigma)$ with $r\varphi(\sigma)$, there
is no reason for the new function $M(\sigma)$ to satisfy ($\varphi
1$).

Let us briefly comment conditions $(\varphi 1)$ and $(\varphi 2)$.  It
is easy to see that they are satisfied when $\varphi(\sigma)=\sigma$,
namely by analytic functions.  In this case
Theorem~\ref{thm:nishihara} provides an alternative proof of the
global existence result for analytic initial data under more
restrictive assumptions on the nonlinearity $m(\sigma)$ (the classical
result in the analytic case only requires $m(\sigma)$ to be continuous
and nonnegative).  More important, assumptions $(\varphi 1)$ and
$(\varphi 2)$ are satisfied when
$\varphi(\sigma)=\sigma/\log(1+\sigma)$, in which case the
corresponding space contains non-analytic functions.  Finally,
assumption $(\varphi 2)$ is not satisfied when
$\varphi(\sigma)=\sigma^{1/s}$ with $s>1$, which means that Gevrey
spaces are never contained in Nishihara's spaces.

\subsection{Our result}\label{sec:result}

In this paper we extend Nishihara's result by replacing assumptions
($\varphi 1$) and ($\varphi 2$) of Theorem~\ref{thm:nishihara} with
the unique assumption (\ref{defn:qa}). Our main result is the following.

\begin{thm}\label{thm:main}
	Let $H$ be a separable Hilbert space, and let $A$ be a nonnegative
	self-adjoint (unbounded) operator on $H$ with dense domain.  Let
	$m:[0,+\infty)\to(0,+\infty)$ be a locally Lipschitz continuous 
	function satisfying the nondegeneracy condition (\ref{hp:ndg}).
	
	Let $\varphi:[0,+\infty)\to[0,+\infty)$ be a continuous and
	strictly increasing function satisfying (\ref{defn:qa}).  Let us
	finally assume that
	\begin{equation}
		(u_{0},u_{1})\in\G_{\varphi,r_{0},3/4}(A)\times
		\G_{\varphi,r_{0},1/4}(A)
		\label{hp:main-data}
	\end{equation}
	for some $r_{0}>0$.
	
	Then problem (\ref{pbm:h-eq}), (\ref{pbm:h-data}) admits a unique
	global solution
	\begin{equation}
		u\in C^{1}\left([0,+\infty);\G_{\varphi,r_{0},3/4}(A)\right)\cap
		C^{0}\left([0,+\infty);\G_{\varphi,r_{0},1/4}(A)\right).
		\label{th:main}
	\end{equation}
\end{thm}

Let us comment our assumptions on the weight $\varphi$, on the initial
data, and on the nonlinearity.

\begin{rmk}
	\begin{em}
		There do exist strictly increasing continuous functions
		$\varphi(\sigma)$ satisfying (\ref{defn:qa}) but not $(\varphi
		1)$.  A nontrivial example is provided in
		section~\ref{sec:spectral-gap}.  A careful inspection of that
		example reveals that not only the function $\sigma\to
		M(\sqrt{\sigma})$ is not convex, but also its convex envelope
		is a constant function (due to the fact that
		$M(\sqrt{\sigma_{k}})=\sqrt{\sigma_{k}}$ on a sequence
		$\sigma_{k}\to +\infty$).
		
		This shows that Theorem~\ref{thm:main}
		is a real extension of Theorem~\ref{thm:nishihara}, and cannot
		be deduced from Theorem~\ref{thm:nishihara} applied with a
		smaller weight which satisfies ($\varphi 1$) and ($\varphi 2$)
		and generates the same functional space.
	\end{em}
\end{rmk}

\begin{rmk}
	\begin{em}
		The ``Sobolev-type'' indices 3/4 and 1/4 of (\ref{th:main})
		are quite usual in the theory of Kirchhoff equations (see for
		example most of the results stated in~\cite{gg:london}).  On
		the other hand, if $\varphi$ grows fast enough (for example if
		$\varphi(\sigma)\geq\log^{2}\sigma$ for every $\sigma\geq 1$),
		then the inclusion $\G_{\varphi,r_{0},\alpha}(A)\subseteq
		\G_{\varphi,r_{1},\beta}(A)$ holds true for every
		$0<r_{1}<r_{0}$ and every $0\leq\alpha\leq\beta$.
		
		We also point out that we allow any $r_{0}>0$ in
		(\ref{hp:main-data}).  This is just because
		condition~(\ref{defn:qa}) doesn't change if we replace
		$\varphi(\sigma)$ with $r_{0}\varphi(\sigma)$.
	\end{em}
\end{rmk}

\begin{rmk}
	\begin{em}
		Concerning the nonlinearity $m(\sigma)$, there is no hope to
		relax assumption (\ref{hp:ndg}) to $m(\sigma)\geq 0$, or the
		Lipschitz continuity assumption to mere continuity.  The
		reason is that some examples presented in~\cite{gg:local} show
		that under these weaker assumptions the Cauchy problem
		(\ref{pbm:h-eq}), (\ref{pbm:h-data}) is not even
		\emph{locally} well posed in classes of quasi-analytic
		functions.
	\end{em}
\end{rmk}

\setcounter{equation}{0}
\section{Proofs}\label{sec:proofs}

\subsection{Technical preliminaries}

In this section we collect some estimates which are crucial in the
proof of our main result.  First of all we remark that assumption
(\ref{defn:qa}) implies in particular that $\varphi(\sigma)$ is
unbounded.  Since it is also continuous and strictly increasing, it
easily follows that $\varphi$, thought as a function
$\varphi:[0,+\infty)\to[\varphi(0),+\infty)$, is invertible.  From now
on we can therefore consider its inverse function
$\varphi^{-1}:[\varphi(0),+\infty)\to[0,+\infty)$.

In the first result we show that assumption (\ref{defn:qa}) implies an
integrability condition on $\varphi^{-1}(\sigma)$ similar to
Nishihara's assumption ($\varphi 2$).

\begin{lemma}\label{lemma:int}
	Let $\varphi:[0,+\infty)\to[0,+\infty)$ be a strictly increasing
	continuous function satisfying (\ref{defn:qa}).  Let
	$\varphi^{-1}:[\varphi(0),+\infty)\to[0,+\infty)$ be its inverse
	function.
	
	Then for every $a>0$, $b\geq 0$, $c>\varphi(0)$ we have that
	\begin{equation}
		\int_{c}^{+\infty}\frac{1}{a\varphi^{-1}(y)+b}\,dy=+\infty.
		\label{th:int-div}
	\end{equation}
\end{lemma}

\subparagraph{\textmd{\emph{Proof}}} 

Assumption (\ref{defn:qa}) is equivalent to say that
$$\int_{d}^{+\infty}\frac{\varphi(y)}{(ay+b)^{2}}\,dy=+\infty$$
for every $a>0$, $b\geq 0$, $d>0$. Let us consider the functions
$$F(x):=\int_{c}^{x}\frac{\varphi(y)}{(ay+b)^{2}}\,dy,
\hspace{3em}
G(x):=\int_{c}^{x}\frac{1}{a\varphi^{-1}(y)+b}\,dy,$$
$$H(x):=G(\varphi(x))-aF(x)-\frac{\varphi(x)}{ax+b},$$
defined for every $x\geq c> \varphi(0)$.  We claim that $H(x)$ is
constant.  If we prove this claim, then (\ref{th:int-div}) easily
follows because 
$$\int_{c}^{+\infty}\frac{1}{a\varphi^{-1}(y)+b}\,dy\
=\ \lim_{x\to +\infty}G(x)\ =\ \lim_{x\to +\infty}G(\varphi(x))\ \geq\
H(c)+a\lim_{x\to +\infty}F(x)\ =$$
$$ =\ H(c)+a\int_{c}^{+\infty}\frac{\varphi(y)}{(ay+b)^{2}}\,dy\ 
=\ +\infty.$$

In order to prove the claim, let us assume first that $\varphi$ is of 
class $C^{1}$. In this case an elementary computation shows that
\begin{equation}
	H'(x)=\frac{\varphi'(x)}{ax+b}-\frac{a\varphi(x)}{(ax+b)^{2}}-
	\frac{\varphi'(x)}{ax+b}+\frac{a\varphi(x)}{(ax+b)^{2}}=0.
	\label{eq:H'}
\end{equation}

If $\varphi$ is not of class $C^{1}$ (and not even absolutely
continuous), then there are at least two standard ways to obtain the
same conclusion.  The first one is to approximate $\varphi(x)$ with a
sequence of strictly increasing functions of class $C^{1}$ and then
passing to the limit.  The second one is recalling that $\varphi$ lies
in $BV_{loc}((0,+\infty))$.  Since $G$ is Lipschitz continuous one can
therefore apply the chain rule in $BV$ and obtain (\ref{eq:H'}) as an
equality between measures instead of functions.\qed

The second result is quite classical.  Roughly speaking, it says that
a solution of the differential inequality (\ref{hp:sub-sol}) cannot
blow up in finite time when the integrability condition
(\ref{hp:int-g}) is satisfied.

\begin{lemma}\label{lemma:ode}
	Let $L_{1}$ be a real number, let
	$g:[L_{1},+\infty)\to(0,+\infty)$ be a positive continuous
	function, and let $y_{0}\geq L_{1}$ be such that
	\begin{equation}
		\int_{y_{0}}^{+\infty}\frac{1}{g(y)}\,dy=+\infty.
		\label{hp:int-g}
	\end{equation}
	
	Let $T>0$, and let $y:[0,T)\to[L_{1},+\infty)$ be a function of class
	$C^{1}$ such that $y(0)=y_{0}$, and
	\begin{equation}
		y'(t)\leq g(y(t))
		\quad\quad
		\forall t\in[0,T).
		\label{hp:sub-sol}
	\end{equation}
	
	Then
	\begin{equation}
		\limsup_{t\to T^{-}}y(t)<+\infty.
		\label{th:limsup}
	\end{equation}
\end{lemma}

\subparagraph{\textmd{\emph{Proof}}} 

Let us consider the function $\Gamma:[L_{1},+\infty)\to\re$ defined by
$$\Gamma(x):=\int_{y_{0}}^{x}\frac{1}{g(y)}\,dy.$$

Assumption (\ref{hp:sub-sol}) is equivalent to say that 
$\left[\Gamma(y(t))\right]'\leq 1$, hence
$$\Gamma(y(t))\leq\Gamma(y(0))+t=t
\quad\quad
\forall t\in[0,T).$$

Due to assumption (\ref{hp:int-g}) we have that $\Gamma$, thought as a
function $\Gamma:[L_{1},+\infty)\to[\Gamma(L_{1}),+\infty)$, is strictly
increasing and invertible, hence 
$$y(t)\leq\Gamma^{-1}(t) \quad\quad
\forall t\in[0,T).$$

At this point (\ref{th:limsup}) easily follows.\qed

The last result is the technical core of this paper. The rough idea 
is that one can estimate an intermediate norm (in this case the sum 
of $a_{k}\lk$) by means of a lower order norm (the sum of $a_{k}$), 
and a higher order norm (the sum of $a_{k}e^{\varphi(\lk)}$). Usually 
such estimates follow from Jensen type inequalities, hence they do
require convexity assumptions as in Nishihara's paper.

Here we prove a result of this type without using convexity.  The
resulting estimates are weaker than the corresponding ones of the
convex case.  Nevertheless they are enough to deduce the a priori
estimates needed in the sequel, and the proof is surprisingly simple.

\begin{prop}\label{lemma:est}
	Let $\varphi:[0,+\infty)\to[0,+\infty)$ be a strictly increasing
	continuous function satisfying (\ref{defn:qa}).  Let
	$\varphi^{-1}:[\varphi(0),+\infty)\to[0,+\infty)$ be its inverse
	function.
	
	Let $\{a_{k}\}$ and $\{\lk\}$ be two sequences of nonnegative real 
	numbers such that
	$$0<E:=\sum_{k=0}^{\infty}a_{k}<+\infty,
	\quad\quad
	F:=\sum_{k=0}^{\infty}a_{k}\max\{\lk,1\}e^{\varphi(\lk)}<+\infty.$$
	
	Then
	\begin{equation}
		\sum_{k=0}^{\infty}a_{k}\lk\leq E\left\{1+\varphi^{-1}
		\left(\varphi(0)+\log\frac{F}{E}\right)\right\}.
	\label{th:basic-ineq}
	\end{equation}
\end{prop}

\subparagraph{\textmd{\emph{Proof}}} 

First of all we remark that $F\geq E>0$, hence the right-hand side of 
(\ref{th:basic-ineq}) is well defined. Let us set for simplicity
$$\alpha:=\varphi^{-1} \left(\varphi(0)+\log\frac{F}{E}\right),$$
and let
$$A:=\left\{k\in\n:\lk<\alpha\right\},\quad\quad\quad
B:=\left\{k\in\n:\lk\geq\alpha\right\}.$$

Let us write
$$\sum_{k=0}^{\infty}a_{k}\lk=\sum_{k\in A}a_{k}\lk+
\sum_{k\in B}a_{k}\lk,$$
and let us estimate the two sums separately.  In $A$ we have that
\begin{equation}
	\sum_{k\in A}a_{k}\lk\leq\alpha\sum_{k\in A}a_{k}\leq\alpha E.
	\label{ineq:A}
\end{equation}

For every $k\in B$ we have that $\varphi(\lk)\geq\varphi(\alpha)=
\varphi(0)+\log(F/E)$, hence
$$e^{\varphi(\lk)}\geq\exp\left(\varphi(0)+\log\frac{F}{E}\right)
\geq\frac{F}{E},$$
and therefore
\begin{equation}
	\sum_{k\in B}a_{k}\lk\leq \frac{E}{F}\sum_{k\in B}a_{k}\lk
	e^{\varphi(\lk)}
	\leq\frac{E}{F}\sum_{k=0}^{\infty}a_{k}\max\{\lk,1\}e^{\varphi(\lk)}
	=E.
	\label{ineq:B}
\end{equation}

Summing (\ref{ineq:A}) and (\ref{ineq:B}) we obtain
(\ref{th:basic-ineq}).  \qed

\subsection{Proof of the main result}

The strategy of the proof is standard for Kirchhoff equations.  First
of all we know that a local solution exists due to classical results.
Then we estimate first order energies using the conserved Hamiltonian.
Finally we prove an a priori estimate on a \emph{higher order} 
energy. This is the key point where Proposition~\ref{lemma:est} plays 
its role. The a priori estimate excludes blow up, and this is enough 
to deduce global existence.

Throughout the proof we assume, without loss of generality, that
$r_{0}=1$.  Indeed the parameter $r_{0}$ in the definition of
$\G_{\varphi,r_{0},\alpha}(A)$ can always be included in
$\varphi(\sigma)$ without changing the fundamental assumption
(\ref{defn:qa}), as previously remarked.

We also assume that 
\begin{equation}
	|u_{1}|^{2}+|A^{1/2}u_{0}|^{2}\neq 0 
	\label{hp:data-neq-0}
\end{equation}
because otherwise the solution is the constant function $u(t)\equiv
u_{0}$, which is clearly globally defined.

Finally, we assume that $m$ is of class $C^{1}$. Indeed, when $m$ is 
just locally Lipschitz continuous, we can approximate it with a 
sequence of smooth functions and then pass all estimates to the 
limit. We spare the reader from the details of this standard argument.

We also point out that the solution is trivially unique because
$m(\sigma)$ is assumed to be Lipschitz continuous (for uniqueness
issues the interested reader is referred to~\cite{gg:uniqueness}).

\paragraph{\textmd{\emph{Maximal local solutions}}}

Due to (\ref{hp:main-data}) we have in particular that
$(u_{0},u_{1})\in D(A^{3/4})\times D(A^{1/4})$.  Therefore the
classical local existence theory (see~\cite{ap,gg:london}) implies
that problem (\ref{pbm:h-eq}), (\ref{pbm:h-data}) admits a unique
local solution
\begin{equation}
	u\in C^{1}\left([0,T);D(A^{1/4})\right)\cap
	C^{0}\left([0,T);D(A^{3/4})\right).
	\label{reg-u-min}
\end{equation}

Moreover, if $[0,T)$ is the maximal interval where this solution is
defined, then either $T=+\infty$, or
\begin{equation}
	\limsup_{t\to T^{-}}|A^{1/4}u'(t)|^{2}+|A^{3/4}u(t)|^{2}=+\infty.
	\label{eq:limsup}
\end{equation}

So we have only to exclude that (\ref{eq:limsup}) holds true.

\paragraph{\textmd{\emph{Standard energy estimates}}}

Let $u$ be any solution of (\ref{pbm:h-eq}), (\ref{pbm:h-data}), with 
regularity prescribed by (\ref{reg-u-min}). Let $u_{k}(t)$ denote the 
components of $u(t)$ with respect to the orthonormal system $e_{k}$ 
(see the simplifying assumptions stated at the beginning  of 
section~\ref{sec:functional}). Let us set
\begin{equation}
	c(t):=m\left(|A^{1/2}u(t)|^{2}\right),
	\label{defn:c}
\end{equation}
and let 
$$E_{k}(t):=|u_{k}'(t)|^{2}+c(t)\lk^{2}|u_{k}(t)|^{2},$$
$$E(t):=|u'(t)|^{2}+c(t)|A^{1/2}u(t)|^{2}=\sum_{k=0}^{\infty}E_{k}(t).$$

We claim that there exist positive constants $L_{1}$ and $L_{2}$ 
such that
\begin{equation}
	L_{1}\leq E(t)\leq L_{2}
	\quad\quad
	\forall t\in[0,T).
	\label{est:et}
\end{equation}

To this end, we consider the usual Hamiltonian 
$$\mathcal{H}(t):=|u'(t)|^{2}+M\left(|A^{1/2}u(t)|^{2}\right),$$
where
$$M(\sigma):=\int_{0}^{\sigma}m(s)\,ds
\quad\quad
\forall\sigma\geq 0.$$

It is well known that $\mathcal{H}(t)$ is constant. By (\ref{hp:ndg}) 
we have that $M(\sigma)\geq\mu\sigma$ for every $\sigma\geq 0$, hence
\begin{equation}
	|A^{1/2}u(t)|^{2}\leq\mu^{-1}\mathcal{H}(t)=\mu^{-1}\mathcal{H}(0),
	\label{est:auu}
\end{equation}
and therefore
$$\mu\leq c(t)\leq\max\left\{m(\sigma):
0\leq\sigma\leq\mu^{-1}\mathcal{H}(0)\right\}=:c_{1}.$$

It follows that
$$E(t)\leq |u'(t)|^{2}+c_{1}|A^{1/2}u(t)|^{2}\leq
c_{2}\mathcal{H}(t)=c_{2}\mathcal{H}(0)=:L_{2}.$$

Similarly, since $M(|A^{1/2}u(t)|^{2})\leq c_{1}|A^{1/2}u(t)|^{2}$, we
have also that 
$$E(t)\geq |u'(t)|^{2}+\mu|A^{1/2}u(t)|^{2}\geq
c_{3}\mathcal{H}(t)=c_{3}\mathcal{H}(0)=:L_{1},$$
where $L_{1}$ is positive due to (\ref{hp:data-neq-0}).  This
completes the proof of (\ref{est:et}).

\paragraph{\textmd{\emph{Fundamental a priori estimate}}}

Let us set
\begin{equation}
	F(t):=\sum_{k=0}^{\infty}E_{k}(t)\max\{1,\lk\}e^{\varphi(\lk)}.
	\label{defn:F}
\end{equation}

We claim that $F(t)$ is well defined for every $t\in[0,T)$, and
\begin{equation}
	\limsup_{t\to T^{-}}F(t)<+\infty.
	\label{limsup-ft}
\end{equation}

To this end, let us first estimate the derivative of $c(t)$. By 
(\ref{defn:c}) and (\ref{est:auu}) we have that
\begin{eqnarray*}
	|c'(t)| & = & \left|m'\left(|A^{1/2}u(t)|^{2}\right)\right|
	\cdot 2\left|\langle A^{1/4}u'(t),A^{3/4}u(t)\rangle\right|  \\
	 & \leq & \max\left\{|m'(\sigma)|:0\leq\sigma\leq
	 \mu^{-1}\mathcal{H}(0)\right\}\cdot
	 \left(|A^{1/4}u'(t)|^{2}+|A^{3/4}u(t)|^{2}\right)\\
	 & \leq & c_{4}\sum_{k=0}^{\infty}\lk E_{k}(t),
\end{eqnarray*}
where the last series converges to a continuous function because we 
already know that $u$ is at least as regular as prescribed by 
(\ref{reg-u-min}). Now we have that
\begin{equation}
	E_{k}'(t)=c'(t)\lk^{2}|u_{k}(t)|^{2}\leq
	\frac{|c'(t)|}{\mu}c(t)\lk^{2}|u_{k}(t)|^{2}
	\leq c_{5}E_{k}(t)\sum_{k=0}^{\infty}\lk E_{k}(t),
	\label{est:ek'}
\end{equation}
hence
\begin{equation}
	E_{k}(t)\leq E_{k}(0)\exp\left(c_{5}
	\int_{0}^{t}\sum_{k=0}^{\infty}\lk E_{k}(\tau)d\tau\right).
	\label{est:ekt}
\end{equation}

On the other hand, from assumption (\ref{hp:main-data}) it is easy to 
deduce that
$$\sum_{k=0}^{\infty}E_{k}(0)\max\{1,\lk\}e^{\varphi(\lk)}<+\infty.$$

Combining with (\ref{est:ekt}) we obtain that the series in 
(\ref{defn:F}) converges, which proves that $F(t)$ is well defined.
Moreover, from (\ref{est:ek'}) we deduce also that $F$ is of class
$C^{1}$, and its derivative satisfies 
$$F'(t)= \sum_{k=0}^{\infty}E_{k}'(t)\max\{1,\lk\}e^{\varphi(\lk)}
\leq c_{5}\left(\sum_{k=0}^{\infty}\lk E_{k}(t)\right)\cdot
F(t).$$

Now we apply Proposition~\ref{lemma:est} with $a_{k}=E_{k}(t)$.  We
obtain that
$$\sum_{k=0}^{\infty}\lk E_{k}(t)\leq E(t)\cdot
\left\{1+\varphi^{-1}\left(\varphi(0)+
\log\frac{F(t)}{E(t)}\right)\right\},$$
hence by (\ref{est:et})
$$\sum_{k=0}^{\infty}\lk E_{k}(t)\leq L_{2}
\left\{1+\varphi^{-1}\left(\varphi(0)+
\log\frac{F(t)}{L_{1}}\right)\right\},$$
and in particular
$$F'(t)\leq c_{6}F(t)
\left\{1+\varphi^{-1}\left(\varphi(0)+
\log\frac{F(t)}{L_{1}}\right)\right\}.$$

Since $F(t)>0$ for every $t\in[0,T)$, from Lemma~\ref{lemma:ode} it
follows that $F(t)$ satisfies (\ref{limsup-ft}) provided that the
function
$$g(y):=c_{6}y\left\{1+\varphi^{-1}\left(\varphi(0)+
\log\frac{y}{L_{1}}\right)\right\},$$
defined for every $y\geq L_{1}$,
satisfies (\ref{hp:int-g}). With the variable change 
$z=\varphi(0)+\log(y/L_{1})$ we obtain that
$$\int_{F(0)}^{+\infty}\frac{1}{g(y)}\,dy=
\int_{c_{7}}^{+\infty}\frac{1}{c_{6}(\varphi^{-1}(z)+1)}\,dz,$$
and the last integral is equal to $+\infty$ due to 
Lemma~\ref{lemma:int} and our assumption (\ref{defn:qa}).

\paragraph{\textmd{\emph{Conclusion}}}

For every $t\in[0,T)$ we have that
$$|A^{1/4}u'(t)|^{2}+|A^{3/4}u(t)|^{2} \leq
c_{8}\sum_{k=0}^{\infty}\lk E_{k}(t) \leq c_{8}F(t).$$
	  
Therefore (\ref{eq:limsup}) cannot hold true because of 
(\ref{limsup-ft}). This is enough to conclude that the solution is 
global. From the estimate on $F(t)$ it follows also that $u$ is as 
regular as required in~(\ref{th:main}).
\qed

\setcounter{equation}{0} 
\section{Connection with spectral-gap solutions}
\label{sec:spectral-gap}

In this section we present some weird results which can be obtained 
using weights $\varphi(\sigma)$ without convexity assumptions. The 
main idea is that a function $\varphi$ can satisfy assumption 
$(\ref{defn:qa})$ even if its growth is very slow (for example 
logarithmic) on a suitable sequence diverging to $+\infty$.

Let $\lk$ be a sequence of nonnegative real numbers such that
$\lambda_{k+1}\geq e^{\lambda_{k}}$ for every $k\in\n$.  Let
$\widetilde{\varphi}:[0,+\infty)\to[0,+\infty)$ be a piecewise
constant function such that
$$\widetilde{\varphi}(x)=\log\lambda_{k+1}
\quad\quad
\forall x\in(\lambda_{k},\lambda_{k+1}]$$
for every $k\in\n$.  Since $\widetilde{\varphi}(x)\geq\lk$ for every
$x\in(\lambda_{k},\lambda_{k+1}]$ and every $k\geq 1$, we have that
\begin{equation}
	\int_{1}^{+\infty}\frac{\widetilde{\varphi}(\sigma)}{\sigma^{2}}\,d\sigma
	\geq\sum_{k=1}^{\infty}\lk\int_{\lk}^{\lambda_{k+1}}
	\frac{d\sigma}{\sigma^{2}}=\sum_{k=1}^{\infty}\left(
	1-\frac{\lk}{\lambda_{k+1}}\right)=+\infty,
	\label{est:series}
\end{equation}
where the last inequality follows from the fact that 
$\lambda_{k+1}\geq 2\lk$ for every $k\geq 1$.

Therefore it is quite simple to modify $\widetilde{\varphi}(\sigma)$ in 
order to obtain a function $\varphi(\sigma)$ which is continuous (or 
even more regular), strictly increasing, satisfies 
$\widetilde{\varphi}(\lk)=\varphi(\lk)$ for all $k\in\n$, and still 
fulfils assumption (\ref{defn:qa}).

Let us assume now that the spectrum of an operator $A$ coincides with
the sequence $\lk$ we have just considered.  Since
$\varphi(\lk)=\log\lk$ for every $k\geq 1$, it follows that
$$\G_{\varphi,r,\alpha}(A)=D(A^{\alpha+r/4}) \quad\quad
\forall\alpha\geq 0,\ \forall r>0.$$

As a consequence, for this operator Theorem~\ref{thm:main} is 
actually a global existence result in Sobolev-type spaces. We have 
thus proved the following result.

\begin{thm}\label{thm:sg-op}
	Let $H$ be a separable Hilbert space, and let $A$ be a nonnegative
	self-adjoint operator on $H$ with dense domain.  Let us assume
	that the spectrum of $A$ consists of a sequence $\lk^{2}$ of
	eigenvalues such that $\lambda_{k+1}\geq e^{\lk}$ for every
	$k\in\n$.
	
	Let $m:[0,+\infty)\to(0,+\infty)$ be a locally Lipschitz
	continuous function satisfying the nondegeneracy condition
	(\ref{hp:ndg}).
	
	Then for every $\alpha>1/4$, and every pair of initial conditions
	$(u_{0},u_{1})\in D(A^{\alpha+1/2})\times D(A^{\alpha})$, problem
	(\ref{pbm:h-eq}), (\ref{pbm:h-data}) admits a (unique) global
	solution
	\begin{equation}
		u\in C^{1}\left([0,+\infty);D(A^{\alpha})\right)\cap
		C^{0}\left([0,+\infty);D(A^{\alpha+1/2})\right).
		\label{th:reg-alpha}
	\end{equation}
\end{thm}

The above result shows that there do exist special operators for which
the Kirchhoff equation is well posed in Sobolev-type spaces.  These
operators are characterized by a fast growing sequence of eigenvalues,
and for this reason we could call them ``spectral-gap operators''.

The same result can also be seen from a different point of view.  The
operator is now more general, but initial data have nonzero components
only with respect to a sequence of special eigenvectors whose
eigenvalues grow fast enough.  We obtain the following result.

\begin{thm}\label{thm:sg-data}
	Let $A$ be a self-adjoint linear operator on a Hilbert space $H$.
	Let us assume that there exist a countable (not necessarily
	complete) orthonormal system $\{e_{k}\}$ in $H$, and a
	sequence $\{\lambda_{k}\}$ of nonnegative real numbers such
	that $\lambda_{k+1}\geq e^{\lk}$ and $Ae_{k}=\lambda_{k}^{2}e_{k}$
	for every $k\in\n$. Let $S\subseteq H$ be the closure of the 
	subspace generated by $\{e_{k}\}$.
	
	Let $m:[0,+\infty)\to(0,+\infty)$ be a locally Lipschitz
	continuous function satisfying the nondegeneracy condition
	(\ref{hp:ndg}).
	
	Let $\alpha>1/4$, and let $(u_{0},u_{1})\in
	D(A^{\alpha+1/2})\times D(A^{\alpha})$ with $(u_{0},u_{1})\in
	S\times S$.
	
	Then problem (\ref{pbm:h-eq}), (\ref{pbm:h-data}) admits a
	(unique) global solution satisfying  (\ref{th:reg-alpha}).
\end{thm}

The proof trivially follows from Theorem~\ref{thm:sg-op} applied in 
$S$. In both cases we could save something on the growth of $\lk$ by 
just asking that
$$\sum_{k=1}^{\infty}\left(\frac{1}{\lk}-
\frac{1}{\lambda_{k+1}}\right)\log\lambda_{k+1}=+\infty,$$
which is what is really required in (\ref{est:series}).

Theorem~\ref{thm:sg-data} reminded us of the recent global existence
results for spectral-gap initial data.  The connection actually exists
because it is possible to show that in most cases elements of $S$ lie
in the spaces introduced by R.\ Manfrin in~\cite{manfrin1,manfrin2}.
On the other hand, the proof given in~\cite{manfrin2} is based on
completely different techniques and requires more regularity both of
initial data ($\alpha\geq 1/2$ instead of $\alpha>1/4$), and of the
nonlinear term ($m(\sigma)$ is assumed to be of class $C^{2}$).  Such
technical restrictions have been recently removed in~\cite{gg:global},
but at the expenses of a faster growth of the sequence of eigenvalues.

In other words Theorem~\ref{thm:sg-op} and Theorem~\ref{thm:sg-data},
in the form they are stated here, don't follow from the theories
developed in~\cite{manfrin1,manfrin2} and~\cite{gg:global}.

\label{NumeroPagine}

\end{document}